\documentclass[final,3p,times,sort]{elsarticle}
\usepackage{amssymb}
\usepackage{amsmath}
\usepackage{graphicx}
\usepackage{graphics}
\usepackage{algorithm}
\usepackage{algorithmic}
\usepackage{expl3}
\usepackage{enumitem}

\usepackage[usenames,dvipsnames,svgnames,table]{xcolor}
\usepackage[pagebackref=true,colorlinks,linkcolor=blue,citecolor=magenta]{hyperref}

\newtheorem{theorem}{Theorem}%[section]
\newtheorem{corollary}[theorem]{Corollary}
\newtheorem{lemma}[theorem]{Lemma}
\newtheorem{proposition}[theorem]{Proposition}
\newtheorem{observation}[theorem]{Observation}
\newtheorem{conjecture}[theorem]{Conjecture}
\newtheorem{remark}[theorem]{Remark}
\newtheorem{preexample}{{\bf Example}}

\newenvironment{example}[1]{\begin{preexample}{\hspace{-0.5
                  em}{\bf }}{\rm #1}}{\end{preexample}}

\newtheorem{preproof}{{\bf Proof}}

\newenvironment{proof}[1]{\begin{preproof}{\rm
              #1}\hfill{$\blacksquare$}}{\end{preproof}}

\newtheorem{presproof}{{\bf Sketch of Proof}}

\newenvironment{sproof}[1]{\begin{presproof}{\rm
              #1}\hfill{$\blacksquare$}}{\end{presproof}}

\newtheorem{claim}{Claim}
\def\F{\mathcal{F}}

\def\ss{\boldsymbol{s}}

\def\ss{\boldsymbol{s}}

\def\cprime{$'$}

\def\zero{\boldsymbol{0}}

\def\ds{\displaystyle}
\def\KG{\operatorname{KG}}

\def\alt{\operatorname{alt}}
\def\HH{\mathcal{H}}

\def\sd{\operatorname{sd}}
\def\mod{\operatorname{mod}}

\def\cd{\operatorname{cd}}
\def\ecd{\operatorname{ecd}}

\journal{Journal}

\begin{document}

\begin{frontmatter}



\title{A new lower bound for the chromatic number of \\
general Kneser hypergraphs}
\author{Roya Abyazi Sani} \ead{roya.abyazisani@shahroodut.ac.ir}
\author{Meysam Alishahi \corref{cor1}}\cortext[cor1]{Corresponding author.}
\ead{meysam$\_$alishahi@shahroodut.ac.ir}
\address{Faculty of Mathematical Sciences, Shahrood University of Technology, Shahrood, Iran.}
\begin{abstract}
A general Kneser hypergraph $\KG^r(\HH)$ is an $r$-uniform hypergraph that somehow encodes the edge intersections of a ground hypergraph $\HH$.
The colorability defect of $\HH$ is a combinatorial parameter providing a lower bound for the 
chromatic number of $\KG^r(\HH)$, which is addressed in a series of works by 
Dol'nikov~[{\it Sibirskii Matematicheskii Zhurnal, 1988}], 
K\v{r}\'{\i}\v{z}~[{\it Transaction of the American Mathematical Society, 1992}], 
and Ziegler~[{\it Inventiones Mathematicae, 2002}].  
In this paper, we define a new combinatorial parameter, the equitable colorability defect of hypergraphs, 
which provides some common improvements upon these works. Roughly speaking,  we propose a new lower bound for 
the chromatic number of general Kneser hypergraphs
which substantially improves  Ziegler's lower bound.  
It is always as good as  Ziegler's lower bound and several 
families of hypergraphs for which the difference between these two lower bounds is arbitrary large are provided. 
This specializes to a substantial improvement 
of the Dol'nikov-K\v{r}\'{\i}\v{z} lower bound for the chromatic number of general Kneser hypergraphs as well. 
Furthermore, we prove a result ensuring the existence 
of a colorful subhypergraph in any proper coloring of general Kneser hypergraphs that strengthens Meunier's result
[{\it The Electronic Journal of Combinatorics,} 2014]. 
\end{abstract}

\end{frontmatter}

% \linenumbers

\section{Introduction and Main Results} 
For two positive integers $n$ and $k$, the two symbols $[n]$ and ${[n]\choose k}$ respectively 
stand for the set $\{1,\ldots, n\}$ and the family of all $k$-subsets of $[n]$.  
A hypergraph $\HH$ is a pair $(V(\HH),E(\HH))$ where $V(\HH)$ -~the {\it vertex set}~- is a finite set and $E(\HH)$ -~the {\it edge set}~- is a family of nonempty subsets of $V(\HH)$. 
For an integer $r\geq 2$, the hypergraph $\HH$ is called {\it $r$-uniform} if each of its edges has the cardinality $r$. 
For a set $N\subseteq V(\HH)$, the {\it induced subhypergraph $\HH[N]$} is the hypergraph whose vertex set is $N$ and whose edge set is $\{e\in E(\HH)\colon e\subseteq N\}$. 

Let $n,k,$ and $r$ be positive integers, where $r\geq 2$ and $n\geq rk$.
The {\it Kneser hypergraph $\KG^r(n,k)$} is an $r$-uniform hypergraph with vertex set 
${[n]\choose k}$ and edge set consisting of all
$r$-tuples of pairwise disjoint members of ${[n]\choose k}$, i.e.,  
$$E(\KG^r(n,k))=\left\{\{e_1,\ldots,e_r\}\colon |e_i|=k, e_i\subseteq [n],\mbox{ and } e_i\cap e_j=\varnothing
\mbox{ for each $i\neq j\in [r]$} \right\}.$$ 
In a break-through, Lov\'{a}sz~\cite{Lovasz1978} determined the chromatic number of 
Kneser graphs $\KG^2(n,k)$ and solved a long-standing conjecture posed by 
Kneser~\cite{Kneser1955}. His proof gave birth to an area of combinatorics which nowadays 
is known as topological combinatorics. Mainly, this area of combinatorics is focused 
on the study of coloring properties of graphs and hypergraphs by using algebraic topological tools. Alon, Frankl, and 
Lov{\'a}sz~\cite{Alonetal1986} extended Lov\'{a}sz's result to 
Kneser hypergraphs $\KG^r(n,k)$ by proving  
$$\chi(\KG^r(n,k))=\left\lceil{n-r(k-1)\over r-1}\right\rceil.$$ 
This result  gave an affirmative answer to a conjecture posed 
by Erd{\H{o}}s~\cite{Erdos1976} as well. Above-mentioned results were generalized in many ways. 
One of the most promising generalizations is the one found by 
Dol'nikov~\cite{Dolnikov1988}, extended by K\v{r}\'{\i}\v{z}~\cite{Kriz, Kriz2000}, and generalized by Ziegler~\cite{Ziegler2002,Ziegler2006}. To state these 
results, we first need to introduce some preliminary notation and definitions. 

Let $\HH$ be a hypergraph and $r$ an integer  where $r\geq 2$.  
The {\it general Kneser hypergraph $\KG^r(\mathcal{H})$} is a 
hypergraph with vertex set $E(\mathcal{H})$ and edge set 
$$E(\KG^r(\mathcal{H})) = \left\{\{e_1,...,e_r\}\colon e_i \in E(\mathcal{H}) \text{ and } 
e_i \cap e_j = \varnothing \text{ for each } i\neq j \in [r] \right\}.$$
For $r=2$, we would use $\KG(\HH)$ rather than $\KG^2(\HH)$. 
Note that the Kneser hypergraph $\KG^r(n,k)$ can be obtained in this way by setting 
$\HH=\left([n],{[n]\choose k}\right)$.  
The {\it $r$-colorability defect of $\HH$}, denoted by $\cd^r(\HH)$,  
is the minimum number of vertices that must be removed from $\HH$ so that the induced 
subhypergraph on the remaining vertices is $r$-colorable. 
Dol'nikov~\cite{Dolnikov1988} (for $r=2$) and K\v{r}\'{\i}\v{z}~\cite{Kriz, Kriz2000}
proved that 
$$\chi\left(\KG^r(\HH)\right)\geq \left\lceil{\cd^r(\HH)\over r-1}\right\rceil.$$
If $\HH=\left([n],{[n]\choose k}\right)$, then this result implies 
$\KG^r(n,k)\geq \left\lceil{n-r(k-1)\over r-1}\right\rceil$
concluding the Alon-Frankl-Lov{\'a}sz result~\cite{Alonetal1986} (one can color $\KG^r(n,k)$ by $\left\lceil{n-r(k-1)\over r-1}\right\rceil$ colors). 

The {\it equitable $r$-colorability defect of $\HH$,} denoted by $\ecd^r(\HH)$, is the minimum 
number of vertices that must be removed from $\HH$ so that the induced subhypergraph on the 
remaining vertices has an equitable $r$-coloring.  Recall that a hypergraph has an
{\it equitable $r$-coloring} if it admits a proper $r$-coloring such that the sizes of its color classes 
differ by at most one, see~\cite{Walte-1973,Yuster-2003}. As an immediate consequence of one of our main results  
(Theorem~\ref{colorfultheorem}), we have the following improvement to  the Dol'nikov-K\v{r}\'{\i}\v{z} lower bound.
\begin{theorem}\label{corollary} 
For any hypergraph $\HH$ and any integer  $r$ with $r\geq 2$, we have 
$$\chi(\KG^r(\HH))\geq \left\lceil \dfrac{\ecd^r(\HH)}{r - 1} \right\rceil.$$
\end{theorem}
The following example illustrates the previous theorem is a 
true improvement to  the Dol'nikov-K\v{r}\'{\i}\v{z} lower bound. However, in the last section, we provide more examples in which the difference between 
the preceding lower bound and the Dol'nikov-K\v{r}\'{\i}\v{z} lower bound can be arbitrary large. \\
  
\begin{example}{
Let $n,k$ and $a$ be integers, where $n> a\geq rk-1$.   
Define $\HH=\HH(n,k,a)$ to be a hypergraph with vertex set $[n]$ and whose edge set is defined as follows: 
$$E(\HH)=\left\{e\subseteq [n]\colon |e|=k\mbox{ and } e\not\subseteq \{n-a+1,\ldots,n\} \right\}.$$ 
}\end{example}\label{exm}
In the last section of this paper, we shall prove that 
$\chi\left(\KG^r(\HH)\right)=\left\lceil \dfrac{n-a}{r-1}\right\rceil.$ 
In fact, we will see that the preceding result simply follows from the lower bound 
stated in Theorem~\ref{corollary}. 
What is interesting about Example~\ref{exm} is that we cannot obtain  
the appropriate lower bound by using the colorability defect of 
$\HH$ (the Dol'nikov-K\v{r}\'{\i}\v{z} lower bound) or even the alternation number of $\HH$. 
It should be noticed that, for any hypergraph $\F$, there is a lower bound for the chromatic 
number of $\KG^r(\F)$ based on the alternation number of $\F$ which surpasses the 
Dol'nikov-K\v{r}\'{\i}\v{z} lower bound 
(definitions and more details will be provided in Section~\ref{comparing}). 
Some other examples comparing these lower bounds will be exhibited in Section~\ref{comparing}.

Let $\HH=\left([n],E(\HH)\right)$ be a hypergraph, $r$ an integer, and $\ss=(s_1,\ldots,s_n)$ 
an integer vector, where $r\geq 2$ and $1\leq s_i< r$ for each $i\in[n]$. 
An $r$-(multi-)set $\left\{N_1,\ldots,N_r\right\}$ is called {\it $\ss$-disjoint} if each 
$i\in[n]$ appears in at most $s_i$ of the sets $N_j$ (we count the repetitions). 
Note that being $(1,\ldots,1)$-disjoint is the same thing as being pairwise disjoint. 
Ziegler~\cite{Ziegler2002} introduced the $\ss$-disjoint general Kneser hypergraph 
$\KG^r_{\ss}(\mathcal{H})$ and the $\ss$-disjoint $r$-colorability defect $\cd^r_{\ss}(\HH)$ 
as generalizations of the general Kneser hypergraph $\KG^r(\HH)$ and the $r$-colorability 
defect $\cd^r(\HH)$, respectively. The {\it $\ss$-disjoint general Kneser hypergraph 
$\KG^r_{\ss}(\mathcal{H})$ } is an $r$-uniform (multi-)hypergraph with vertex set 
$E(\HH)$ and edge set 
$$E(\KG^r_{\ss}(\HH)) = 
\Big\{\{e_1,\ldots,e_r\}\colon e_1,\ldots ,e_r\in E(\HH) \text{ and } \{e_1,\ldots ,e_r\} 
\text{ is $\ss$-disjoint }\Big\}.$$ 
Depending on $\ss$, an edge of 
$\KG^r_{\ss}(\HH)$ is not necessarily a set and it might be a multi-set of size $r$. 
Similar to the Kneser hypergraph $\KG^r(n,k)$, the {\it $\ss$-disjoint Kneser hypergraph} $\KG^r_{\ss}(n,k)$ is defined to be 
$\KG^r_{\ss}\left([n],{[n]\choose k}\right)$. 
Also, the 
{\it $\ss$-disjoint $r$-colorability defect of $\HH$,} denoted by $\cd^r_{\ss}(\HH)$, is defined as follows: 
$$\cd_{\ss}^r(\HH) =\bar{n}-\max\left\{\ds\sum_{j=1}^r|N_j|\colon N_j\subseteq [n],\;  
\{N_1,\ldots,N_r\}\mbox{ is $\ss$-disjoint and $E(\HH[N_j])=\varnothing$ for each } j\in [r]\right\},$$
where $\bar{n}=\bar{n}(\ss)=\sum\limits_{i=1}^n s_i$.

Ziegler~\cite{Ziegler2002,Ziegler2006} generalized  the Dol'nikov-K\v{r}\'{\i}\v{z} lower bound to 
$\ss$-disjoint general Kneser hypergraphs. Indeed, he proved that 
%the following inequality holds 
\begin{equation}\label{Ziegler's inequality}
\chi(\KG^r_{\ss}(\HH))\geq \left\lceil{\cd_{\ss}^r(\HH)\over r-1}\right\rceil,
\end{equation}
provided that $\max \limits _{i \in [n]}  s_i <\mu(r)$, where $\mu(r)$ is the largest prime factor 
of $r$. Note that Ziegler's lower bound immediately implies the Dol'nikov-K\v{r}\'{\i}\v{z} 
lower bound since for $\ss=(1,\ldots,1)$, we have $\cd^r_{\ss}(\HH)=\cd^r(\HH)$ and 
$\KG^r_{\ss}(\HH) = \KG^r(\HH)$. 
Ziegler, using his lower bound, determined the chromatic number of the $\ss$-disjoint 
general Kneser hypergraph $\KG^r_{\ss}(n,k)$ for some integer vectors 
$\ss$.  

An $r$-uniform hypergraph $\HH$ is called {\it $r$-partite}  
if its vertex set can be partitioned into subsets (parts) $V_1,\ldots,V_r$  so that each of its edges 
intersects each part $V_i$ in exactly one vertex. An $r$-uniform $r$-partite hypergraph is said to 
be {\it complete} if it contains all possible edges. For the case $r=2$, such a graph is called a 
{\it complete bipartite graph}. 

Generalizing Lov{\'a}sz's result, several results concern the 
existence of large colorful complete bipartite subgraphs in properly colored graphs 
when some condition related to a topological lower bound on the chromatic number is satisfied, 
e.g.~\cite{Alishahi&Hajiabolhassan2015, Changetal2013, Chen2011, 
Simonyietal2013, Simonyi&Tardos2006, Simonyi&Tardos2007}.   
In this regard, Simonyi and Tardos~\cite{Simonyi&Tardos2007} proved that for an arbitrary 
hypergraph $\HH$, any proper coloring of the general Kneser graph $\KG(\HH)$ contains a 
multicolored complete bipartite subgraph $K_{\lfloor{t\over 2}\rfloor,\lceil{t\over 2}\rceil}$ of 
order $t=\cd^2(\HH)$ in which the colors alternate on two parts of this bipartite subgraph with respect to their 
natural order.
The existence of colorful subhypergraphs in hypergraphs coloring was first studied  
by Meunier~\cite{Meunier2014}. Meunier generalized the Simonyi-Tardos result to general 
Kneser hypergraphs $\KG^p(\HH)$ provided that $p$ is prime. 
Actually, he proved that for a prime number $p$ and  
a hypergraph $\HH$, any proper coloring of $\KG^p(\HH)$ contains a $p$-uniform $p$-partite 
subhypergraph  with parts $V_1,\ldots,V_p$ satisfying the following properties:  
\begin{itemize}
\item $\ds\sum_{i=1}^p|V_i|=\cd^p(\HH)$,
\item the values of $|V_i|$'s differ by at most one, i.e., $||V_i|-|V_j||\leq 1$ for each $i,j\in[p]$, and 
\item for any $j\in[p]$, the vertices in $V_j$ get distinct colors. 
\end{itemize} 
The aforementioned result by Simonyi and Tardos ensures the existence  of a colorful subgraph with $\cd^2(\HH)$ vertices 
in any proper coloring of $\KG(\HH)$.  
Some extensions of this result to more general topological settings can be found in~\cite{Alishahi&Hajiabolhassan2015,Simonyietal2013}.  
These results as well as Meunier's were extended to the case of uniform hypergraphs in~\cite{Alishahi2017}. 

Let $r$ and $n$ be two integers and $\ss=(s_1,\ldots,s_n)$ an integer vector, where $1\leq s_i<r$. For the sets $N_1,\ldots,N_r\subseteq [n]$, the (multi-)set 
$\{N_1,\ldots,N_r\}$ is called {\it equitable }  if the values of $|N_j|$'s differ by at most one.  
Also, it is called  {\it equitable  $\ss$-disjoint} if it is simultaneously equitable and $\ss$-disjoint. 
For a hypergraph $\HH=([n], E(\HH))$, define the {\it equitable  $\ss$-disjoint 
$r$-colorability defect of $\HH$}, denoted by $\ecd_{\ss}^r(\HH)$, to be the following quantity: 
 \begin{eqnarray*}
 \ecd^r_{\ss}(\HH) =\bar{n}-\max\left\{\ds\sum_{j=1}^r|N_j|\colon N_j\subseteq [n],  \{N_1,\ldots,N_r\}\mbox{ is equitable $\ss$-disjoint and $E(\HH[N_j])=\varnothing$ for each } j\in [r]\right\}.
\end{eqnarray*}
Note, here, $\bar{n}= \sum \limits _{i = 1} ^ n s_i$. Throughout the paper, for $\ss=(1,\ldots,1)$,  since $\ecd^r(\HH)$ and $\ecd_{\ss}^r(\HH)$ are 
the same, we use  $\ecd^r(\HH)$ rather than $\ecd_{\ss}^r(\HH)$.

The next theorem is the first main result of the paper which not only extends 
Meunier's result~\cite{Meunier2014} to general $\ss$-disjoint  Kneser hypergraphs 
$\KG^p_{\ss}(\mathcal{H})$ but also improves it to equitable $r$-colorability defect. 
This theorem also results in a lower bound for the chromatic number of 
$\KG^r_{\ss}(\mathcal{H})$ based on the equitable $\ss$-disjoint $r$-colorability defect 
$\ecd^r_{\ss}(\mathcal{H})$ which surpasses Ziegler's lower bound~\cite{Ziegler2002,Ziegler2006}.

\begin{theorem} \label{maintheorem}
Let $\HH=([n],E(\HH))$ be a hypergraph, $p$ a prime number, and $\ss=(s_1,\ldots,s_n)$  
a positive integer vector, where $1\leq s_i< p$ for each $i\in [n]$. Any proper coloring of the 
$\ss$-disjoint general Kneser hypergraph $\KG^p_{\ss}(\mathcal{H})$ contains some 
subhypergraph whose vertex set can be partitioned into parts $V_1,\ldots,V_p$ satisfying the following properties;
\begin{itemize}
\item $\sum\limits_{i=1}^p|V_i|=\ecd^p_{\ss}(\HH)$,
\item $\{e_1,\ldots,e_p\}$ forms an edge of $\KG^p_{\ss}(\mathcal{H})$ for each choice of $e_i\in V_i$, 
\item the values of $|V_i|$'s differ by at most one, i.e., $||V_i|-|V_j||\leq 1$ for each $i,j\in[p]$, and 
\item for any $j\in[p]$, the vertices in $V_j$ get distinct colors. 
\end{itemize}
\end{theorem}
If $\ss=(1,\ldots,1)$, then Theorem~\ref{maintheorem} implies Meunier's theorem
in a stronger form (using $\ecd^r(\HH)$ instead of $\cd^r(\HH)$). Also, since any color 
appears in at most $p-1$ vertices of each edge of $\KG^p_{\ss}(\mathcal{H})$, 
this theorem results in $\chi(\KG^p_{\ss}(\HH)) \geq 
\left\lceil \dfrac{\ecd^p_{\ss}(\HH)}{p - 1} \right \rceil$ provided that $p$ is prime. 
The following theorem, applying this result, provides an improvement of Ziegler's theorem~\cite{Ziegler2002,Ziegler2006}. 
\begin{theorem}\label{colorfultheorem}
Let $\HH=([n],E(\HH))$ be a hypergraph, $r$ an integer, $\ss=(s_1,\ldots,s_n)$  
a positive integer vector, and  $\mu(r)$ be the largest prime factor of $r$. If  $1\leq s_i< \mu(r)$ for each $i\in [n]$,  then   
$$\chi(\KG^r_{\ss}(\HH))\geq \left\lceil \dfrac{\ecd^r_{\ss}(\HH)}{r - 1} \right\rceil.$$
\end{theorem}
Since we always have $\ecd^r_{\ss}(\HH)\geq \cd^r_{\ss}(\HH)$, 
setting $\ss=(1,\ldots,1)$ leads to the improvement of the 
Dol'nikov-K\v{r}\'{\i}\v{z} lower bound stated in Theorem~\ref{corollary}. \\

\noindent{\bf Plan.} 
The rest of this paper is organized as follows. In Section~\ref{proof}, we first introduce some topological  
tools which will be needed elsewhere in the paper.  Next, we prove 
Theorem~\ref{maintheorem} and deduce Theorem~\ref{colorfultheorem} from 
Theorem~\ref{maintheorem} by reducing Theorem~\ref{colorfultheorem} to the case of prime $r$. 
We conclude the paper in Section~\ref{compare} by 
a discussion on comparing the equitable colorability defect with some other combinatorial  
parameters providing lower bounds for the chromatic number of general 
Kneser hypergraphs.  Indeed, we will build some families of hypergraphs $\HH$ in which the difference between 
the lower bound introduced in Theorem~\ref{colorfultheorem} and some other well-known lower bounds for the chromatic number of general Kneser hypergraphs $\KG^r(\HH)$
can be arbitrary large.

\section{Proofs of the main results}\label{proof}
This section is devoted to the proofs of Theorems~\ref{maintheorem} and~\ref{colorfultheorem}. 
First, we review some basic definitions and tools. Next, we prove Theorem~\ref{maintheorem}. 
To prove Theorem~\ref{colorfultheorem}, we reduce this theorem 
to the case of prime $r$, which has been already proved by the discussion after 
Theorem~\ref{maintheorem}. 

\subsection{Basic Tools}\label{basictools}
Here, a brief review of some notation and definitions which will be used throughout this 
section is provided though it is assumed that the reader has basic knowledge on topological 
combinatorics (for more, see~\cite{matousek2003}). 

For an integer $r$, by $ Z_r$, we refer to the cyclic multiplicative group of order $r$ and with 
generator $\omega$, i.e., $ Z_r=\{\omega,\omega^2,\ldots,\omega^r\}$. 
A {\it simplicial complex}, considered as a combinatorial object or a topological space, 
is a pair $(V, K)$, where $V$ - the {\it vertex set} - is a finite set and  $K$ - 
the {\it simplex set} - is a hereditary system of nonempty subsets of $V$, that is, 
if $F \in K$ and $\varnothing\neq F' \subseteq F$, then $F' \in K$. 
By a simplicial complex $K$, we mean the simplicial complex $(V, K)$, where 
$V = \bigcup \limits_{F \in K} F$. Each set in $K$ is called a {\it simplex} of $K$.
The {\it dimension of $K$}, denoted by $\dim(K)$, equals $\max \limits _{F \in K} |F| -1$. 
The {\it first barycentric subdivision of the simplicial complex $K$}, denoted by $\sd K$, is  
the order-complex obtained from the poset consisting of all simplices in $K$ 
ordered by inclusion. The join of two simplicial complexes $C$ and $K$, denoted by $C*K$,  
is a simplicial complex with vertex set $V(C) \uplus V (K)$ and simplex set  
$$\Big\{F\subseteq V(C) \uplus V (K)\colon F\neq\varnothing, F\cap V(C)\in  C\cup\{\varnothing\} \text{ and } F\cap V(K) \in K\cup\{\varnothing\}\Big\}.$$ 
%$$\Big\{F_1 \uplus F_2 : F_1 \in C\cup\{\varnothing\} \text{ and } F_2 \in K\cup\{\varnothing\}\Big\}\setminus\big\{\varnothing\big\}.$$}
In view of the definition, for three simplicial complexes $K,L$, and $M$, the two simplicial complexes $K*(L*M)$ and $(K*L)*M$ are the same simplicial complexes (up to a natural relabeling of their vertices). This allows us to use $K*L*M$ for both of $K*(L*M)$ and $(K*L)*M$ whenever we do not care about the names of the vertices.  
The join of $n$ disjoint copies of $K$ is denoted by $K^{*n}$. By renaming the vertices of $K^{*n}$, we can see  
$K^{*n}$ as a simplicial complex with vertex set $V(K)\times [n]$ such that for each vertex $(v,i)\in V(K)\times [n]$, the index $i$ indicates that $v$ is considered as a vertex of the 
$i$-th copy of $K$. 
Let $p$ be a prime and $s$ be an integer such that $1\leq s<p$. The simplicial complex $\sigma^{p-1}_{s-1}$ 
has $Z_p$ as the vertex set and its simplex set consists of all nonempty subsets of 
$Z_p$ with size at most $s$. 

For an integer vector $\ss=(s_1,\ldots,s_n)\in [p-1]^n$,
the {\it $\ss$-disjoint $p$-fold join $(\sigma^{n-1})^p_{\ss}$} is a simplicial complex whose vertex set is $Z_p\times [n]$ and its simplex set consists of all nonempty subsets $A$ of $Z_p\times [n]$ such that for each $i\in [n]$, the pair $(\varepsilon,i)$ is present in $A$ for at most $s_i$ different $\varepsilon\in Z_p$, i.e.,
$$\left|\left\{\varepsilon\in Z_p\colon (\varepsilon, i)\in A\right\}\right|\leq s_i.$$ One can simply check that $(\sigma^{n-1})^p_{\ss}=\sigma^{p-1}_{s_1-1}*\cdots*\sigma^{p-1}_{s_n-1}$. 

Let $C$ and $K$ be two simplicial complexes.
By a {\it simplicial map} $f : C \rightarrow K$, we mean a map from $V(C)$ to $V(K)$ 
such that the image of each simplex of $C$ is a simplex of $K$. 
The simplicial complex $C$ is a {\it simplicial $ Z_p$-complex} if 
$ Z_p$ acts on it and moreover, 
the map $\varepsilon: C \rightarrow C$ which maps $v$ to $\varepsilon\cdot v$  %$v \mapsto \varepsilon\cdot v$ 
is a simplicial map for each $\varepsilon \in  Z_p$. 
A simplicial $ Z_p$-complex is called {\it free} if there is no fixed simplex under the simplicial map 
made by each $\varepsilon \in  Z_p \setminus \{\boldsymbol{1}\}$,
where $\boldsymbol{1}$ is the identity element of the group $ Z_p$.
For two simplicial $ Z_p$-complexes $C$ and $K$,  a simplicial map $f : C \rightarrow K$ is called 
{\it $ Z_p$-equivariant} if $f(\varepsilon\cdot v) = \varepsilon\cdot f(v)$ for each 
$\varepsilon \in  Z_p$ and $v \in V (C)$. Simply, any $ Z_p$-equivariant simplicial map is called 
{\it simplicial $ Z_p$-map}. Dold's theorem~\cite{Dold1982} is a transformation group extension of the Borsuk-Ulam theorem which, in particular, asserts that if there is a 
simplicial $Z_p$-map from a free $Z_p$ space $C$ to a free $Z_p$ space $K$, then the dimension of $K$ is larger than the connectivity of $C$.
For the definition of connectivity of topological spaces, we refer the reader to the book by  Matou{\v{s}}ek~\cite{matousek2003}.  
It is known that $\sigma^{p-1}_{s-1}$ is a free  
$(s-2)$-connected simplicial complex (with the natural $Z_p$-action). This along with the connectivity lemma for joins which asserts that the connectivity of the join of two topological spaces is at least the summation of their connectivity plus $2$ concludes that $(\sigma^{n-1})^p_{\ss}$ (and its barycentric subdivision)  
has the connectivity $\bar{n}-2=\sum\limits_{i=1}^ns_i -2$ (see the topological proof of Lemma~5.3 in~\cite{Ziegler2002}). 

\subsection{Proof of Theorem \ref{maintheorem} }\label{theorem for prime}
In our approach, we utilize the function $l(-)$ introduced 
by~\cite{Alishahi2017} and some sign functions 
proposed by Meunier~\cite{Meunier2014}. 
Although, our technique is somehow similar to that used in~\cite{Alishahi2017}, we will develop it to work here.
To be more specific, the desired conclusion follows by applying Dold's theorem to a simplicial $Z_p$-map defined via a given proper coloring of the general Kneser hypergraph
between two simplicial $ Z_p$-complexes.  However, in~\cite{Alishahi2017}, we did~not care about the equitability  and moreover we did~not deal with 
multi-hypergraphs. Indeed, rather than in~\cite{Alishahi2017}, we here work with different simplicial complexes, use sign functions differently, and moreover, 
we define the simplicial $Z_p$-map in a different way. 
Before starting the proof, we give some preliminaries. 
Let $q$ be a positive integer and  $p$ a prime number. 
For a simplex $\tau\in\left(\sigma^{p-1}_{p-2}\right)^{*q}$, 
define $h(\tau) = \min \limits_{\varepsilon \in Z_p} |\tau ^\varepsilon|,$ 
where $\tau ^\varepsilon =\big\{(\varepsilon, j) : (\varepsilon, j)\in \tau\big\}$.   
%Let $\tau$ be a simplex in $\left(\sigma^{p-1}_{p-2}\right)^{*n}$.  Set  $h(\tau) = \min \limits_{\varepsilon \in Z_p} |\tau ^\varepsilon|.$ Also, for each $\varepsilon \in Z_p$, define $\tau ^\varepsilon =\{(\varepsilon, j) : (\varepsilon, j)\in \tau\}$.  
Also, set $$l(\tau) = p\cdot h(\tau) + |\{\varepsilon \in Z_p : |\tau ^ \varepsilon| > h(\tau)\}|.$$  
It should be emphasized that we shall use the two functions $l(-)$ and $h(-)$ several times 
during the proof.

\begin{proof}{
Let $\HH=([n],E(\HH))$ be a hypergraph, $p$ a prime number, and $\ss=(s_1,\ldots,s_n)$  
a positive integer vector, where $1\leq s_i< p$ for each $i\in [n]$. Also, let 
$c:E(\HH)\longrightarrow [t]$ be a proper coloring of $\KG_{\ss}^p(\HH)$. 
When $\ecd_{\ss}^p(\HH)=0$, the statement of Theorem~\ref{maintheorem} is clearly true. 
For the rest of the proof, we shall assume that $\ecd_{\ss}^p(\HH)>0$. 

Hereafter, for simplicity of notation, we set 
$K=\sigma_{s_1-1}^{p-1}*\cdots*\sigma_{s_n-1}^{p-1}$. 
For a positive integer $a \in [n]$, let $W_a$ denote the set consisting of all simplices 
$\tau \in  K$ such that 
$|\tau ^\varepsilon| \in \{0, a\}$ for each $\varepsilon \in Z_p$. 
Also, for a positive integer $b \in [t]$, let $U_b$ be the set consisting of all simplices 
$\tau \in \left(\sigma^{p-1}_{p-2}\right)^{*t}$ such that $|\tau ^\varepsilon| \in \{0, b\}$ for each 
$\varepsilon \in Z_p$.
Define $W = \bigcup\limits _{a = 1}^n W_a$ and $U=\bigcup\limits _{b = 1}^t U_b$. 
Choose three arbitrary $Z_p$-equivariant maps $s_0 : \sigma^{p-1}_{p-2} \rightarrow Z_p$, $s_1 : W \rightarrow Z_p$, and $s_2 : U \rightarrow Z_p$ 
by choosing one representative in each orbit (note that this is possible because 
$ Z_p$ acts freely on each of $\sigma^{p-1}_{p-2}$, $W$, and $U$).

For each simplex $\tau$ of $K$ and for each $\varepsilon\in Z_p$,  define 
$$\tau_\varepsilon=\left\{j\in[n]\colon (\varepsilon,j)\in \tau\right\}.$$
Note that $\tau_\varepsilon\subseteq [n]=V(\HH).$
Further, define 
$$l(\HH) = \max \left\{l(\tau)\colon \tau \in K\text{ such that } E(\HH[\tau _\varepsilon]) = \varnothing \mbox{ for each }  
\varepsilon \in Z_p \right\}.$$
When $\{\tau \in K\colon  E(\HH[\tau _\varepsilon]) = \varnothing \mbox{ for each }  \varepsilon \in Z_p\}=\varnothing$, i.e., each singleton is an edge of $\HH$, 
we define $l(\HH)=0$.  
One should notice that $\bar{n}-l(\HH)=\ecd_{\ss}^r(\HH)$.  
Now, set $$m=l(\HH)+\max\big\{l(\tau_c)\colon \tau\in K\mbox{ and } l(\tau)>l(\HH)\big\},$$ 
where, for each simplex $\tau$ in $K$, 
$$\tau_c=\big\{(\varepsilon,c(e))\in  Z_p\times [t]\colon e\in E(\HH)\mbox{ and }e\subseteq 
\tau_\varepsilon \big\}.$$
Since $\ecd_{\ss}^r(\HH)=\bar{n}-l(\HH)>0$, the set $\{l(\tau_c)\colon \tau\in K\mbox{ and } l(\tau)>l(\HH)\big\}$
is~not empty and consequently $m$ is well defined. 
In view of the definition of $K=(\sigma^{n-1})^p_{\ss}$, for each simplex $\tau$ of $K$, 
the set $\{\tau_\varepsilon\colon \varepsilon\in  Z_p\}$ is $\ss$-disjoint. 
This observation along with the fact that $c$ is a proper coloring of $\KG^p_{\ss}(\HH)$ 
implies $\tau_c$ is either empty or a simplex of $(\sigma_{p-2}^{p-1})^{*t}$. 
Indeed, if $l(\tau)>l(\HH)$, then 
$\tau_c$ is a simplex of $(\sigma_{p-2}^{p-1})^{*t}$. This observation plays a crucial role in our proof. 

In what follows, we first define a map 
$\gamma\colon \sd K\longrightarrow  Z_p^{*m}.$ 
Let $\tau$ be an arbitrary simplex in $K$. Define $\gamma(\tau)$ as follows.  
\begin{itemize}
\item[{\rm a)}] For $l(\tau)\leq l(\HH)$,    
		      we consider the two following cases. 
                    \begin{enumerate}
                         \item If $h(\tau)=0$, 
                         		 then define 
                         		 $\gamma(\tau)=(s_0(\bar{\tau}),l(\tau))$, where 
                         	         $$\bar{\tau}=\{\varepsilon\in Z_p\colon  
	         		 \tau^\varepsilon=   \varnothing\}\in \sigma^{p-1}_{p-2}.$$
			\item If $h(\tau)> 0$, then define 
				$\gamma(\tau)=\left(s_1(\bar\tau), l(\tau)\right)$, 
				where 
				$$\bar{\tau}=\displaystyle \bigcup_{\{\varepsilon\in Z_p\colon 
				|\tau^\varepsilon|=h(\tau)\}} \tau^\varepsilon\in W.$$
		    \end{enumerate}
\item[{\rm b)}] For $l(\tau)> l(\HH)$, the set $\tau_c$ is not empty and thus it is a simplex of 
		      $(\sigma_{p-2}^{p-1})^{*t}$. Now, we consider two different cases. 
                    \begin{enumerate}
                         \item  If $h(\tau_c)=0$, then define 
                         		 $\gamma(\tau)=(s_0(\bar{\tau}_c),l(\HH)+l(\tau_c))$, where 
                         	         $$\bar{\tau}_c=\{\varepsilon\in Z_p\colon  
	         		 {\tau_c}^\varepsilon=   \varnothing\}\in \sigma^{p-1}_{p-2}.$$
			\item If $h(\tau_c)> 0$, then define 
				$\gamma(\tau)=\left(s_2(\bar{\tau}_c), l(\HH)+l(\tau_c)\right)$, 
				where 
				$$\bar{\tau}_c=\displaystyle \bigcup_{\{\varepsilon\in Z_p\colon  
				|{\tau_c}^\varepsilon|=h(\tau_c)\}} {\tau_c}^\varepsilon\in U.$$
		    \end{enumerate}
\end{itemize} 
\begin{claim}\label{claim1}
The map $\gamma$ is a simplicial $Z_p$-map
from $\sd K$ to $ Z_p^{*m}$.
\end{claim}
\begin{proof}{
First, note that all the functions $s_0(-)$, $s_1(-)$, and $s_2(-)$ are $Z_p$-equivariant
maps. Accordingly, it is clear that  $\gamma$ is a $Z_p$-equivariant map as well. 
For a contradiction, suppose that  $\gamma$ is~not a simplicial map. 
In view of the simplices in $\sd K$ and $ Z_p^{*m}$, we conclude 
there are two simplices $\tau$ and $\tau'$  in $K$
such that $\tau\subsetneq\tau'$,  $\gamma(\tau)=(\varepsilon_1,\beta)$, and 
$\gamma(\tau')=(\varepsilon_2,\beta)$, 
where $\varepsilon_1\neq \varepsilon_2$ and $\beta\in [m]$. 
In view of the definition of $\gamma$, it is clear that 
it is~not possible to have $l(\tau)\leq l(\HH)$ and $l(\tau')> l(\HH)$, simultaneously.  
Therefore, we have either $l(\tau)\leq l(\tau')\leq l(\HH)$ or $l(\tau')\geq l(\tau)>l(\HH)$, 
which will be discussed separately in what follows.
\begin{itemize}
\item[{\rm I)}] $l(\tau)\leq l(\tau')\leq l(\HH).$ 
Clearly, in view of the definition of $\gamma$ in this case, we should have 
$l(\tau)=l(\tau')$. We consider the two following cases.\\
     \begin{enumerate}
		\item If $h(\tau)=h(\tau')=0$, then since $\tau\subsetneq\tau'$ and 
			$$\varepsilon_1=s_0(\{\varepsilon\in Z_p\colon           
			\tau^\varepsilon=   \varnothing\})\neq s_0(\{\varepsilon\in Z_p\colon    
			{\tau'}^\varepsilon=   \varnothing\})=\varepsilon_2,$$  
			we conclude $ \{\varepsilon\in Z_p\colon   
			{\tau'}^\varepsilon=   \varnothing\}\subsetneq\{\varepsilon\in Z_p\colon  
			\tau^\varepsilon=   \varnothing\}$. This implies that 
			$$l(\tau')=p-|\{\varepsilon\in Z_p\colon  
			{\tau'}^\varepsilon=   \varnothing\}|>p-|\{\varepsilon\in Z_p\colon  
			\tau^\varepsilon=   \varnothing\}|=l(\tau),$$ 
			a contradiction.\\
		\item For $h(\tau')>0$, we know that $$ l(\tau)=p\cdot h(\tau)+
			 |\{\varepsilon\in Z_p\colon  |\tau^\varepsilon|>h(\tau)\}|
			 \mbox{ and } l(\tau')=p\cdot h(\tau')+|\{\varepsilon\in Z_p\colon 
			 |{\tau'}^\varepsilon|>h(\tau')\}|.$$  
			 The facts  $l(\tau)=l(\tau')$ and  
			 $$\max\left\{|\{\varepsilon\in Z_p\colon  |\tau^\varepsilon|>h(\tau)\}|, 
			 |\{\varepsilon\in Z_p\colon 
			 |{\tau'}^\varepsilon|>h(\tau')\}|\right\}\leq p-1$$ 
			 imply that  $h=h(\tau)=h(\tau')$ and 
			 $$|\{\varepsilon\in Z_p\colon  |\tau^\varepsilon|>h\}| 
			 = |\{\varepsilon\in Z_p\colon  
			 |{\tau'}^\varepsilon|>h\}|.$$ 
			 In view of  the equations   
			 $$\varepsilon_1=s_1(\displaystyle\bigcup_{\{\varepsilon\in Z_p\colon  
			 |\tau^\varepsilon|=h\}} \tau^\varepsilon)\neq 
			 s_1(\displaystyle\bigcup_{\{\varepsilon\in Z_p\colon  
			 |{\tau'}^\varepsilon|=h\}} {\tau'}^\varepsilon)=\varepsilon_2,$$ 
			 we must have $$\displaystyle\bigcup_{\{\varepsilon\in Z_p\colon  
			 |\tau^\varepsilon|=h\}} \tau^\varepsilon\neq 
			 \displaystyle\bigcup_{\{\varepsilon\in Z_p\colon  
			 |{\tau'}^\varepsilon|=h\}} {\tau'}^\varepsilon.$$ 
			 This inequality in combination with the facts that  
			 $\tau\subsetneq \tau'$ and $h=\ds\min_{\varepsilon\in  Z_p}
			 |\tau^\varepsilon|=\ds\min_{\varepsilon\in  Z_p}|{\tau'}^\varepsilon|$ 
			 imply that  $$\left\{\varepsilon\in Z_p\colon  |{\tau'}^\varepsilon|=h\right\}
			 \subsetneq \left\{\varepsilon\in Z_p\colon |{\tau}^\varepsilon|=h\right\},$$ 
			 which is impossible.
    \end{enumerate}
\item[{\rm II)}] $l(\tau')\geq l(\tau)> l(\HH).$ 
Note both $\tau_c$ and ${\tau'}_c$ are simplices in $(\sigma_{p-2}^{p-1})^{*t}$. 
Clearly, in view of the definition of $\gamma$ in this case, we should have 
$l(\tau_c)=l({\tau'}_c)$.  Also, note that since $\tau\subsetneq \tau'$, we have 
$\tau_c\subseteq {\tau'}_c$. Furthermore, the equality cannot happen, since otherwise,  in view of the definition of $\gamma$, 
we must have $\varepsilon_1=\varepsilon_2$, which contradicts the way $\tau$ and $\tau'$ have been chosen.     
Similar to the previous case, we will deal with the two different 
cases $h(\tau_c)=h({\tau'}_c)=0$ and $h({\tau'}_c)>0$. 
     \begin{enumerate}
		\item If $h(\tau_c)=h({\tau'}_c)=0$, then since $\tau_c\subsetneq\tau'_c$ and 
			$$\varepsilon_1=s_0(\{\varepsilon\in  Z_p\colon           
			\tau_c^\varepsilon=   \varnothing\})\neq s_0(\{\varepsilon\in  Z_p\colon    
			{\tau'}_c^\varepsilon=   \varnothing\})=\varepsilon_2,$$  
			we have    $ \{\varepsilon\in  Z_p\colon   
			{\tau'}_c^\varepsilon=   \varnothing\}\subsetneq\{\varepsilon\in  Z_p\colon  
			\tau_c^\varepsilon=   \varnothing\}$. This implies that 
			$$l({\tau'}_c)=p-|\{\varepsilon\in  Z_p\colon  
			{\tau'}_c^\varepsilon=   \varnothing\}|>p-|\{\varepsilon\in  Z_p\colon  
			\tau_c^\varepsilon=   \varnothing\}|=l(\tau_c),$$ 
			a contradiction.\\
		\item For $h({\tau'}_c)>0$, we know that $$ l(\tau_c)=p\cdot h(\tau_c)+
			 |\{\varepsilon\in Z_p\colon  |\tau_c^\varepsilon|>h(\tau_c)\}|
			 \mbox{ and } l({\tau'}_c)=p\cdot h({\tau'}_c)+|\{\varepsilon\in Z_p\colon  
			 |{\tau'}_c^\varepsilon|>h({\tau'}_c)\}|.$$  
			 The facts  $l(\tau_c)=l({\tau'}_c)$ and  
			 $$\max\left\{|\{\varepsilon\in Z_p\colon  |\tau_c^\varepsilon|>h(\tau_c)\}|, 
			 |\{\varepsilon\in Z_p\colon  
			 |{\tau'}_c^\varepsilon|>h({\tau'}_c)\}|\right\}\leq p-1$$ 
			 conclude  $h=h(\tau_c)=h({\tau'}_c)$ and 
			 $$|\{\varepsilon\in Z_p\colon  |\tau_c^\varepsilon|>h\}| 
			 = |\{\varepsilon\in Z_p\colon  
			 |{\tau'}_c^\varepsilon|>h\}|.$$ 
			 According to  
			 $$\varepsilon_1=s_2(\displaystyle\bigcup_{\{\varepsilon\in Z_p\colon  
			 |\tau_c^\varepsilon|=h\}} \tau_c^\varepsilon)\neq 
			 s_2(\displaystyle\bigcup_{\{\varepsilon\in Z_p\colon  
			 |{\tau'}_c^\varepsilon|=h\}} {\tau'}_c^\varepsilon)=\varepsilon_2,$$ 
			 we must have $$\displaystyle\bigcup_{\{\varepsilon\in Z_p\colon  
			 |\tau_c^\varepsilon|=h\}} \tau_c^\varepsilon\neq 
			 \displaystyle\bigcup_{\{\varepsilon\in Z_p\colon  
			 |{\tau'}_c^\varepsilon|=h\}} {\tau'}_c^\varepsilon.$$ 
			 This inequality and the facts that  
			 $\tau_c\subsetneq {\tau'}_c$ and $h=\ds\min_{\varepsilon\in  Z_p}
			 |\tau_c^\varepsilon|=\ds\min_{\varepsilon\in  Z_p}|
			 {\tau'}_c^\varepsilon|$ 
			 result in $$\left\{\varepsilon\in Z_p\colon |{\tau'}_c^\varepsilon|=h\right\}
			 \subsetneq \left\{\varepsilon\in Z_p\colon |{\tau}_c^\varepsilon|=h\right\},$$ 
			 which is impossible.
    \end{enumerate}
\end{itemize}
}\end{proof}
\begin{claim}
There is a simplex $\tau\in K$  for which $l({\tau}_c)\geq \ecd_{\ss}^p(\HH).$ 
\end{claim}
\begin{proof}{
By Claim~\ref{claim1}, it has already been noted that $\gamma$ is a simplicial $ Z_p$-map from $\sd K$ to 
$ Z_p^{*m},$ where $$m=l(\HH)+\max\left\{l(\tau_c)\colon \tau\in K\mbox{ and } l(\tau)>l(\HH)\right\}.$$
Accordingly, in view of Dold's theorem~\cite{Dold1982}, the dimension of $ Z_p^{*m}$ 
must be strictly larger than the connectivity of $\sd K$; that is 
$m-1>\bar{n}-2$, which implies  
$m\geq \bar{n}$ (see the discussion at the end of Subsection~\ref{basictools}).  Consequently, from the definition of $m$, 
there is a simplex $\tau\in K$  for which $l(\HH)+l({\tau}_c)\geq \bar{n}.$ 
Equivalently, we have 
$$l({\tau}_c)\geq \bar{n}-l(\HH)=\ecd_{\ss}^p(\HH),$$ 
as desired.  
}\end{proof}
Let $\tau$ be a simplex for which $l({\tau}_c)\geq \ecd_{\ss}^p(\HH)$. For simplicity of notation, set $h=h(\tau_c)$. 
For each $i\in [p]$, if $|{\tau_c}^{\omega^i}|=h$, then define  
$C_i$ to be the $h$-set $\{c_{i,1},\ldots,c_{i,h}\}\subseteq [t]$ for which $\{\omega^i\}\times C_i={\tau_c}^{\omega^i}$, 
otherwise, (if $|{\tau_c}^{\omega^i}|>h$) let 
$C_i=\{c_{i,1},\ldots,c_{i,h+1}\}\subseteq [t]$ be an $(h+1)$-set such that 
$\{\omega^i\}\times C_i\subseteq{\tau_c}^{\omega^i}$. 
Now, for each $i\in [p]$, if $|{\tau_c}^{\omega^i}|=h$, then  
define $V_i=\left\{e_{i,1},\ldots,e_{i,h} \right\}$ such that  $e_{i,j}\subseteq \tau_{\omega^i}$ 
and $c(e_{i,j})=c_{i,j}$ for each $j\in[h]$, otherwise, define  
$V_i=\left\{e_{i,1},\ldots,e_{i,h+1} \right\}$ such that  $e_{i,j}\subseteq \tau_{\omega^i}$ 
and $c(e_{i,j})=c_{i,j}$ for each $j\in[h+1]$. 
Clearly, 
$$\sum\limits_{i=1}^p|V_i|=\sum\limits_{i=1}^p|C_i|=l(\tau_c)\geq \ecd^p_{\ss}(\HH).$$
We have already noticed that the set 
$\{\tau_\varepsilon\colon \varepsilon\in Z_p\}$ is $\ss$-disjoint,
therefore, for each choice of $e_i\in V_i$, the set $\{e_1,\ldots,e_p\}$ is $\ss$-disjoint 
as well and consequently $\{e_1,\ldots,e_p\}$ is an edge of $\KG^p_{\ss}(\mathcal{H})$. 
It is then trivial that the subhypergraph 
$\KG^p_{\ss}(\mathcal{H})[\bigcup\limits_{i=1}^pV_i]$ is the desired subhypergraph. 
}\end{proof}

\subsection{Proof of Theorem~\ref{colorfultheorem}}
In this subsection, we reduce Theorem~\ref{colorfultheorem} to the case of prime $r$, which 
is known to be true owing to Theorem~\ref{maintheorem} (see the discussion after 
Theorem~\ref{maintheorem}). 
The idea of this reduction is originally due to 
K\v{r}\'{\i}\v{z}~\cite{ Kriz2000}, which has been used in other papers as well, e.g.~\cite{Alishahi&Hajiabolhassan2015,HaMe16,Ziegler2002,Ziegler2006}.

For a hypergraph $\mathcal{H} = ([n], E(\HH))$ and positive integers $r$ and $C$, define 
$\mathcal{T} _{\mathcal{H}, C, r}$ 
to be a hypergraph with vertex set $[n]$ and edge set 
\begin{eqnarray*}
E(\mathcal{T} _{\mathcal{H}, C, r}) = \ds\Big\{ V \subseteq [n]\colon \ecd^r(\mathcal{H}[V]) > (r-1)C\Big\}.
\end{eqnarray*}
\begin{lemma} \label{auxiliary lemma}
Let $r'$, $r''$ be two positive integers and $\ss=(s_1,\ldots,s_n)$
be an integer vector, where $1\leq s_i<r''$.
For any hypergraph $\mathcal{H}=([n],E(\HH))$, the following inequality holds
\begin{eqnarray}
\ecd^{r'r''}_{\ss} (\mathcal{H}) \leq r''(r' - 1)C + \ecd^{r''}_{\ss}(\mathcal{T}_{\mathcal{H}, C, r'}).
\end{eqnarray}
\end{lemma}
\begin{proof}{
For the ease of use, set $\mathcal{T}=\mathcal{T}_{\mathcal{H}, C, r'}$.
In view of the definition of $\ecd^{r''}_{\ss}(\mathcal{T})$, there exists a family 
of equitable $\ss$-disjoint subsets of $[n]$, say $\{V_1,..., V_{r''}\}$, 
such that $\sum\limits_{i = 1}^{r''} |V_i| = \sum\limits_{i = 1}^{n} s_i - 
\ecd^{r''} _{\ss} (\mathcal{T})$ and $E(\mathcal{T}[V_i]) = \varnothing$ for each $i\in[r'']$.
Clearly, this implies that $V_i\not\in E(\mathcal{T})$ for each $i\in[r'']$ and consequently,  
$\ecd^{r'}(\mathcal{H}[V_i]) \leq (r'-1)C$ for each $i\in[r'']$.
Thus, for each $i\in[r'']$, there is a family $\{V_{i1},...,V_{ir'}\}$ of equitable disjoint 
subsets of $V_i$ such that   
$$\sum \limits^{r'} _{j = 1} |V_{ij}| = |V_i| - \ecd^{r'}
(\mathcal{H}[V_i])\geq |V_i| - (r'-1)C$$ and
$E(\mathcal{H}[V_{ij}]) = \varnothing$ for all $j\in[r']$.
One can simply check that the family  $\left\{V_{ij}\colon 1\leq i\leq r''\; \&\; 1\leq j\leq r'\right\}$   
is an  equitable $\ss$-disjoint family of subsets of $[n]$ and 
moreover,  $E(\mathcal{H}[V_{ij}]) = \varnothing$ for each $ i \in [r'']$ and $j \in [r']$.
Accordingly, 
$$
\begin{array}{ccl}
\ecd^{r'r''}_{\ss}(\mathcal{H}) & \leq & \sum\limits_{i=1}^n s_i- \sum\limits_{i = 1}^{r''}\sum\limits_{j = 1}^{r'} |V_{ij}| \\
					    &   \leq    &  \sum\limits_{i=1}^n s_i -\sum\limits_{i = 1}^{r''}|V_i| + r''(r'-1)C\\
					    &     =     &  \ecd^{r''}_{\ss}(\mathcal{T}) + r''(r'-1)C

\end{array}
$$
which completes the proof. 
}\end{proof}
\begin{lemma}\label{reduction}
Let $r, r', r''$ be positive integers, where $r',r''\geq 2$, $r=r'r''$.  
Also, let $\ss=(s_1,\ldots,s_n)$ be a positive integer vector, 
where $1\leq s_i<r''$ for each $i\in[r'']$.
If Theorem~\ref{colorfultheorem} holds for $r'$ and $\ss'=(1,\ldots,1)$ and also 
for $r''$ and $\ss$, then it holds for $r$ and $\ss$.
\end{lemma}
\begin{proof}{
For contradiction sake, suppose that there is a proper coloring 
$c : E(\mathcal{H}) \rightarrow [C]$ of $\KG^{r}_{\ss}(\mathcal{H})$ for which 
$\ecd ^{r} _{\ss} (\mathcal{H}) > (r-1)C$. 
Applying Lemma \ref{auxiliary lemma} leads us to   
$$(r'r'' - 1) C < \ecd^{r'r''}_{\ss}(\mathcal{H}) \leq r''(r' -1) C +\ecd^{r''}_{\ss}(\mathcal{T}_{\mathcal{H},C,r'}),$$ which immediately implies  
$(r'' - 1)C < \ecd ^{r''}_{\ss}(\mathcal{T}_{\mathcal{H},C,r'}).$
Since Theorem~\ref{colorfultheorem} holds for $r''$ and $\ss$, 
the preceding observation concludes  
$$\chi(\KG^{r''}_{\ss}(\mathcal{T}_{{\mathcal{H},C,r'}}))>C.$$

On the other hand, in view of the definition of $\mathcal{T}_{{\mathcal{H},C,r'}}$, 
for each $e\in E(\mathcal{T}_{\mathcal{H},C,r'})$, we have 
$\ecd ^{r'} (\mathcal{H}[e]) > (r' -1) C$.
Since Theorem \ref{colorfultheorem} holds for $r'$ and $\ss'=(1,\ldots,1)$, this implies 
$\chi (\KG^{r'}(\mathcal{H}[e]))>C$.
Consequently the assignment of colors to the vertices of $\KG^{r'} (\mathcal{H}[e])$ 
through the coloring $c$ cannot be a proper coloring. 
Therefore, for each edge $e \in E(\mathcal{T}_{\mathcal{H},C,r'})$, 
there is at least one monochromatic edge in $\KG^{r'} (\mathcal{H}[e])$. 
Now, for each edge $e \in E(\mathcal{T}_{\mathcal{H},C,r'})$,
set $h(e)$ to be the largest 
color amongst the colors of monochromatic edges in $\KG^{r'}(\mathcal{H}[e])$.
Since $c$ is a proper coloring  of $\KG^r_{\ss} (\mathcal{H})$, the map 
$h: E(\mathcal{H}) \rightarrow C$ is a proper coloring of $\KG^ {r''}_{\ss}
(\mathcal{T}_{\mathcal{H},C,r'})$, which implies  
$\chi(\KG^ {r''}_{\ss}(\mathcal{T}_{\mathcal{H},C,r'}))\leq C,$ a contradiction. 
To see this, contrary to the claim, suppose that $\{f_1,\ldots,f_{r''}\}$ is an edge of  
$\KG^ {r''}_{\ss}(\mathcal{T}_{\mathcal{H},C,r'})$ such that 
$h(f_1)=\cdots=h(f_{r''})=a$. By the definition of $h(-)$, for each $j\in[r'']$,  there is 
a monochromatic edge $\{e_{1,j},\ldots,e_{r',j}\}$ of $\KG^{r'} (\mathcal{H}[f_j])$ 
for which we have $c(e_{1,j})=\cdots=c(e_{r',j})=a$. 
One can simply see that $\left\{e_{i,j}\colon i\in[r']\;\&\; j\in[r'']\right\}$ is a monochromatic edge 
of $\KG^r_{\ss} (\mathcal{H})$, contradicting the fact that $c$ is a proper coloring for
$\KG^r_{\ss} (\mathcal{H})$. 
%Now, since Theorem~\ref{colorfultheorem} holds for  $r''$ and $\ss_1$, we also have 
%$$
%\begin{array}{ccl}
%\chi(\KG^ {r''}_{\ss_1}(\mathcal{T}_{\mathcal{H},C,r',\ss_1})) 
%		& \geq & \dfrac{\ecd^{r''}_{\ss} (\mathcal{T})}{r'' - 1}\\
%		& \geq & \dfrac{\ecd_{\ss}^{r'r''}(\mathcal{H})-r''(r'-1)C}{s(r''-1)}\\
%		&   >    & \dfrac{(r'r''-1)C-r''(r'-1)C}{s(r''-1)}\\
%		&   =    & \dfrac{C}{s},
%\end{array}
%$$which is a contradiction.
}\end{proof} 
Here, we apply Lemma~\ref{reduction} to complete the proof of Theorem~\ref{colorfultheorem}.
First, note that, by the discussion after Theorem~\ref{maintheorem}, we know that  
Theorem~\ref{colorfultheorem} is true for $\ss=(1,\ldots,1)$ and the prime values of $r$. Thus, applying Lemma~\ref{reduction} with $\ss=\ss'=(1,\ldots,1)$, 
Theorem~\ref{colorfultheorem} holds for $\ss=(1,\ldots,1)$ and any positive integer $r\geq 2$. 
%Therefore, for $\ss'=(1,\ldots,1)$, Theorem~\ref{colorfultheorem} holds for any $r\geq 2$. 
Now, in Lemma~\ref{reduction}, we set $r''$ to be the largest prime factor of $r$. Again, in view of the discussion after 
Theorem~\ref{maintheorem}, we know that Theorem~\ref{colorfultheorem} holds for $r''$ and $\ss$. 
Consequently, applying Lemma~\ref{reduction} gives the desired conclusion.

\section{\bf Comparing equitable colorability defect with colorability defect and alternation number}\label{compare} 
In this section, we compare equitable colorability defect of hypergraphs with colorability 
defect and alternation number of them, two other combinatorial parameters providing lower bounds for 
the chromatic number of general Kneser hypergraphs. 
%{\bgreen It should be noticed that even though all examples presented in this section can be simply modified to the general case of $\ss$, for the ease of reading, we prefer to state them in the simpler case $\ss=(1,\ldots,1)$.} 

\subsection{Comparing equitable colorability defect with colorability defect}  
By the definitions of $\cd_{\ss}^r(\HH)$ and $\ecd_{\ss}^r(\HH)$, it is apparent 
$\ecd_{\ss}^r(\HH)\geq \cd_{\ss}^r(\HH)$.
There are several examples ensuring that not only this inequality might be strict but also 
the difference between $\ecd_{\ss}^r(\HH)$ and $\cd_{\ss}^r(\HH)$
%, i.e., $\ecd_{\ss}^r(\HH)- \cd_{\ss}^r(\HH)$, 
can be arbitrary large. To see this, let $m$ and $n$ be positive integers, where $m\geq r(n+1)$. 
Consider an arbitrary hypergraph $\F$  with $n$ vertices and let 
$U=\{u_1,\ldots,u_m\}$ be a set disjoint from $V(\F)$. 
Define $\HH$ to be a hypergraph with  vertex set 
$V(\F)\cup U$ and edge set $E(\HH)=E(\F)\cup \left\{\{u,v\}\colon u\in U\mbox{ and } v\in V(\F)  \right\}.$ 
Since any independent set of $\HH$ is either a subset of $U$ or an independent set of $\F$, 
it is simple to check that $\cd^r(\HH)= \cd^{r-1}(\F)$ and $\ecd^r(\HH)=n$, which results in 
$\ecd^r(\HH)-\cd^r(\HH)=n-\cd^{r-1}(\F)$. 
Now, one can consider several hypergraphs $\F$ for which $n-\cd^{r-1}(\F)$ is 
arbitrary large. By the preceding construction, in particular, we can build some $r$-partite graphs which are also discussed in more detail in the following proposition. 
\begin{proposition}
If we set $\HH$ to be  the complete $r$-partite graph $K_{t,...t,T}$ with $T\geq rt + r -1$, 
then $\cd^r(\HH) = 0$ while $\chi(\KG^r(\HH))={\ecd^r(\HH)\over r-1}=t$. 
\end{proposition}
\begin{proof}{
Let $V(\HH)=V_1\uplus\cdots \uplus V_r$ such that each $V_i$ is independent, $|V_1|=\cdots=|V_{r-1}|=t$, and $|V_r|=T$.  
First note that $\HH$ is $r$-colorable (it is $r$-partite) which clearly implies $\cd^r(\HH) = 0$. 
Let $(S_1,\ldots,S_r)$ be an equitable $r$-coloring of $\HH\big[\bigcup\limits_{i=1}^r S_i\big]$. To prove $\ecd^r(\HH) = (r-1)t$, we need to show that 
$|\bigcup\limits_{i=1}^r S_i|\leq T$ and the equality can be achieved, which will be concluded from the following two cases:  
If there is some $i\in[r]$ such that $S_i\subseteq V_j$ and $j\leq r-1$, then since $(S_1,\ldots,S_r)$ is an equitable partition for  
$\bigcup\limits_{i=1}^r S_i$,  we have $|\bigcup\limits_{i=1}^r S_i|\leq t+(r-1)(t+1)\leq T$. Otherwise, if $S_i\subseteq V_r$ for each $i\in[r]$, then 
$\bigcup\limits_{i=1}^r S_i\subseteq V_r$ which implies $|\bigcup\limits_{i=1}^r S_i|\leq |V_r|=T$. 
It is clear that we can have the equality in latter case. 
In view of Theorem~\ref{corollary}, to finish the proof, it suffices to show that $\KG^r(\HH)$ is $t$-colorable.  
To this end, let $V(\HH)\setminus V_r=U_1\uplus\cdots\uplus U_{t}$ such that $|U_j|=r-1$ for each $j\in[t]$. 
Now, for each edge $e\in E(\HH)$, define $c(e)$ to be the minimum $j$ such that  $e\cap U_j\neq\varnothing.$ 
One can check that  $c$ is a proper $t$-coloring of $\KG^r(\HH)$, as desired. 
}\end{proof}

Let us remind that $\mu(r)$ is the largest prime factor of an integer $r\geq 2$. 
Ziegler~\cite{Ziegler2002,Ziegler2006}, generalizing the Alon-Frankl-Lov{\'a}sz theorem~\cite{Alonetal1986}, determined the chromatic number of some families of $\ss$-disjoint Kneser hypergraphs $\KG_{\ss}^r(n,k)$. For positive integers $n,k,r,$ and $s$ with $k\geq 2$, $ns\geq rk$, and $1\leq s<\mu(r)$, he proved  
$$\chi(\KG_{(s,\ldots,s)}^r(n,k))=\left\lceil{ns-r(k-1)\over r-1}\right\rceil,$$  
provided that $s$ divides $r-1$. Note that if $s=1$, then we have the Alon-Frankl-Lov{\'a}sz theorem.

Let $n,k,r,$ and $a$ be positive integers and $\ss\in[r-1]^n$ a positive integer vector with $k,r\geq 2$, $n>a$, and  $\sum\limits_{i=1}^ns_i\geq rk$. 
Define $\HH(n, k, a)$ to be a hypergraph with vertex set 
$[n]$ and edge set $\left\{ e \in {[n] \choose k}: e \nsubseteq \{n-a+1,\ldots,n\}\right\}$.  
Also, set $\KG^r_{\ss}(n,k,a)=\KG^r_{\ss}(\HH(n, k, a))$. 
Note that $\KG^r_{\ss}(n,k,a)=\KG^r_{\ss}(n,k)$ when $a\leq k-1$.  
In what follows, we compare the equitable colorability defect and colorability defect of $\HH(n,k,a)$ and prove that the difference between these two quantities can be arbitrary large. Moreover, extending the latter result by Ziegler, we study the chromatic number of $\KG_{(s,\ldots,s)}^r(n,k,a)$ in some cases.

\begin{lemma}\label{upper}
Let $n,k,r,$ and $a$ be positive integers and $\ss=(s,\ldots,s)\in[r-1]^n$ a positive integer vector with $k,r\geq 2$, $n> a$, and  $ns\geq rk$. 
Then, we have 
$$\chi(\KG_{\ss}^r(n,k,a))\leq \left\lceil{1\over \lfloor{r-1\over s}\rfloor}{ns- \max\{rk-1,sa+r-1\}\over s}\right\rceil+1.$$
\end{lemma}
\begin{proof}{
The proof is similar to that of Lemma~3.1 in~\cite{Ziegler2002}.  
Set $P=\lfloor{r-1\over s}\rfloor$ and $M=\left\lceil{1\over \lfloor{r-1\over s}\rfloor}{ns- \max\{rk-1,sa+r-1\}\over s}\right\rceil$. 
For each vertex $e\in {[n]\choose k}$, define 
$$c(e)=\min\left\{\left\lceil{\min e\over P} \right\rceil, M+1\right\}.$$
Note that $c$ is a coloring with color set $[M+1]$. To complete the proof, it suffices to show that this coloring is proper. 
Contradictory, suppose that $\{e_1,\ldots,e_r\}$ is an edge of $\KG^r_{\ss}(n,k,a)$ with $c(e_1)=\cdots=c(e_r)=\ell$. 
First note that  $\ell\leq M$ is impossible. Indeed, since $Ps=\lfloor{r-1\over s}\rfloor s\leq r-1$, at most $r-1$ of the edges $e_i$'s can receive a color $\ell\in[M]$. 
If $\ell=M+1$, then $\min e_i> MP$ for each $i\in[r]$ which concludes  
$\bigcup\limits_{i=1}^{r}e_i\subseteq \{MP+1,\ldots, n\}$.  We thus have   
$${rk\over s}\leq |\bigcup\limits_{i=1}^{r}e_i|\leq n-MP\leq {\max\{rk-1,sa+r-1\}\over s},$$
which implies $sa+r-1\geq rk$ and consequently $MP+1\geq n-a-\lfloor{r-1\over s}\rfloor+1$.  
Since $\KG_{\ss}^r(n,k,a)$ has no vertex contained in $\{n-a+1,\ldots,n\}$, each 
$e_i$ must have some element in $\left\{n-a-\lfloor{r-1\over s}\rfloor+1,\ldots, n-a\right\}$. On the other hand, since $\{e_1,\ldots,e_r\}$ is an $(s,\ldots,s)$-disjoint family, each $j\in[n]$ is 
present in at most $s$ of the sets $e_i$ which concludes $s\lfloor{r-1\over s}\rfloor\geq r$, a contradiction. 
}\end{proof}

By the following observation, we will compute the $\ss$-disjoint $r$-colorability defect and the equitable  $\ss$-disjoint 
$r$-colorability defect of $\HH(n,k,a)$ when $\ss=(s,\ldots,s)\in[r-1]^n$.   

\begin{observation}\label{ecrH(n,k,a)}
Let $n,k,r,$ and $a$ be positive integers and $\ss=(s,\ldots,s)\in[r-1]^n$ a positive integer vector with $k,r\geq 2$, $n>a$, and  $ns\geq rk$. 
Then 
$$\begin{array}{l}%\label{ecdcompare}
\cd^r_{\ss}(\HH(n,k,a))  =    
		\begin{cases}  
		                        ns - r(k-1) \quad & a\leq k-1\\ \\
				       \max\big\{0,ns-(r-s)(k-1)-as\big\} \quad & a\geq k
		\end{cases}
\end{array}$$
and 
$$\begin{array}{l}%\label{ecdcompare}
\ecd^r_{\ss}(\HH(n,k,a))  =    
		\begin{cases}  
		                        ns - r(k-1) \quad& a\leq k-1\\ \\
		                        ns - r(k-1)-\left \lfloor\dfrac{as}{k}\right\rfloor \quad & sk\leq as\leq rk-2\\ \\
				       ns - as \quad & as\geq rk-1. 
		\end{cases}
\end{array}$$
\end{observation}
\begin{sproof}{
The proof of Observation~\ref{ecrH(n,k,a)} is simple but technical. To ease the reading, we just sketch the proof of the second equality. 
Let $\{N_1,\ldots,N_r\}$ be an equitable $(s,\ldots,s)$-disjoint family of subsets of $[n]$ such that $\ecd^r_{\ss}(\HH(n,k,a))=ns-\sum\limits_{i=1}^r|N_i|$ and 
$\HH(n,k,a)\big[N_j\big]=\varnothing$ for each $j\in[r]$.  First note that  $\min\limits_{j\in[r]}|N_j|\geq k-1$ and also, for each $N_j$ with $|N_j|\geq k$, 
we must have $N_j\subseteq A=\{n-a+1,\ldots,n\}$. Consequently, one can check that $\sum\limits_{j=1}^r|N_j|$ will be maximized  whenever the number of $N_j$'s with 
$N_j\subseteq A$ and $|N_j|\geq k$ is maximized.  Now, one can simply verify the desired equations. 
}\end{sproof}%\qed\end{proof}
By the next corollary, generalizing Ziegler's result~\cite{Ziegler2002,Ziegler2006}, we compute the chromatic number of $\KG^r_{\ss}(n,k,a)$ in some cases. 
\begin{corollary}\label{corchrom}
Let $n,k,r,$ and $a$ be positive integers and $\ss=(s,\ldots,s)\in[\mu(r)-1]^n$ be a positive integer vector with $k,r\geq 2$, $n>a$, and  $ns\geq rk$. 
If $s$ divides $r-1$, then 
$$\chi(\KG^r_{\ss}(n,k,a)) =   \left\lceil{ns- \max\{r(k-1),as\}\over r-1}\right\rceil$$
provided that either $a\leq k-1$ or $as\geq rk-1$. 
\end{corollary}
\begin{proof}{
Since $s$ divides $r-1$, by Theorem~\ref{colorfultheorem} and Lemma~\ref{upper}, we have 
$$\left\lceil{\ecd^r_{\ss}(\HH(n,k,a))\over r-1}\right\rceil\leq \chi(\KG_{\ss}^r(n,k,a))\leq \left\lceil{ns- \max\{r(k-1),as\}\over r-1}\right\rceil.$$
Now, the proof follows by Observation~\ref{ecrH(n,k,a)}. 
}\end{proof}  
Under the same assumption as in Corollary~\ref{corchrom}, if $s$ divides $r-1$ and $ks\leq as<rk-1$, then, 
by Lemma~\ref{upper}, Theorem~\ref{colorfultheorem}, and Observation~\ref{ecrH(n,k,a)}, we have   
$$\left\lceil{ns - r(k-1)-\left \lfloor\dfrac{as}{k}\right\rfloor\over r-1}\right\rceil\leq \chi(\KG_{\ss}^r(n,k,a))\leq \left\lceil{ns-\max\{r(k-1),as\}\over r-1}\right\rceil.$$
Therefore, using Theorem~\ref{colorfultheorem}, we are~not able to determine the chromatic number of $\KG_{\ss}^r(n,k,a)$ in this case. As an interesting  question, one may ask for the chromatic number of 
$\KG_{\ss}^r(n,k,a)$ provided that $ks\leq as<rk-1$ and $s$ divides $r-1$. In Subsection~\ref{comparing}, we present more evidences supporting the supposition that the chromatic number of $\KG_{\ss}^r(n,k,a)$ is equal to the preceding upper bound provided that $ks\leq as<rk-1$ and $s$ divides $r-1$. 

\subsection{Comparing equitable colorability defect with alternation number} \label{comparing}
Throughout this subsection, assume $\ss=(1,\ldots,1)$.   
Here, we first define the alternation number of hypergraphs. 
%Let $ Z_r  = \big\{\omega ^ j\colon j   \in [r] \big\}$ be a cyclic and multiplicative group of order $r$ with generator $\omega$. 
For an $X=(x_1,...,x_n) \in ( Z_r\cup\{0\})^n$, 
the subsequence $x_{i_1},\ldots,x_{i_m}$ of nonzero coordinates of $X$ is 
called {\it alternating} if any two consecutive terms of this subsequence are different.
In other words, the sequence  $x_{i_1},\ldots,x_{i_m}$ is called alternating whenever  
$1\leq i_1<\cdots<i_m\leq n$, $x_{i_j}\neq 0$ for each $j\in[m]$ and  $x_{i_j}\neq x_{i_{j+1}}$ 
for each $j\in[m-1]$. For each $X\in ( Z_r\cup\{0\})^n$, we define $\alt(X)$ to be the 
longest alternating subsequence of $X$. Also, we define $\alt(\zero)=0$. 
For each $X\in ( Z_r\cup\{0\})^n$ and for each $i\in[r]$, set 
$X^i=\{j\in[n]\colon x_j=\omega^i\}.$ By abuse of notation, we can write $X=(X^1,\ldots,X^r)$.

Let $\HH=(V(\HH),E(\HH))$ be a hypergraph. 
Also, consider a bijection $\sigma:[n]\longrightarrow V(\HH)$. The {\it $r$-alternation number of 
$\HH$ with respect to the bijection $\sigma$}, denoted by $\alt^r_{\sigma}(\HH)$, is the maximum possible $k$ for which there is an 
$X\in ( Z_r\cup\{0\})^n$
with $\alt(X)=k$ such that $E(\HH[\sigma(X^i)])=\varnothing$ for each $i\in[r]$, i.e.,
$$\alt^r_{\sigma}(\HH)=\max\left\{\alt(X)\colon X\in ( Z_r\cup\{0\})^n\mbox{ and }  
E(\HH[\sigma(X^i)])=\varnothing\mbox{ for each } i\in[r]\right\}.$$
Define 
$$\alt^r(\HH)=\min\limits_{\sigma}\alt^r_\sigma(\HH),$$ 
where the minimum is taken over all bijections $\sigma:[n]\longrightarrow V(\HH)$.   
The second present author and Hajiabolhassan~\cite[Theorem~3]{Alishahi&Hajiabolhassan2015}, 
using the $Z_p$-Tucker lemma, proved that 
\begin{eqnarray}\label{alishahi&haji inequality}
\chi(\KG^r(\mathcal{H})) \geq \left\lceil\dfrac{|V(\mathcal{H})| - \alt^r(\mathcal{H})}{r-1}\right\rceil.
\end{eqnarray}
They also computed the chromatic number of several families of hypergraphs using this 
lower bound, see~\cite{Alishahi&Hajiabolhassan2015}.  One can simply see that the 
preceding lower bound surpasses the Dol'nikov-K\v{r}\'{\i}\v{z} lower bound. 

Note that, setting $\ss=(1,\ldots,1)$, we have already computed $\cd^r(\HH(n,k,a))$ and $\ecd^r(\HH(n,k,a))$ by Observation~\ref{ecrH(n,k,a)}. 
However, for computing $\alt^r(\HH(n,k,a))$, we need to put more effort. 
\begin{lemma}\label{altcompare}
Let $n,k,r,$ and $a$ be positive integers such that $k,r\geq 2$ and $n\geq \max\{a+1,rk\}$. Then 
$$\begin{array}{l}
n-\alt^r(\HH(n,k,a))   \leq     
			    \begin{cases} 
			    			n- \max\{2(k-1),\;        a\} \quad & r =2\\ \\
						n- \max\{3(k-1),\;   a+1\}  \quad & r =3\\ \\
                                                  n- \max\{r(k-1),\;    a+b\}  \quad  & r \geq 4,
                      	   \end{cases}
\end{array}$$
where $b=\min\left\{n-a,(r-2)(k-1)\right\}$. 
\end{lemma}
\begin{proof}{ 
For simplicity of notation, set $\HH=\HH(n,k,a)$. 
We will first show that $\alt^r_\sigma(\HH)\geq r(k-1)$.
Then, we will prove that $\alt^2_\sigma(\HH)\geq a$, $\alt^3_\sigma(\HH)\geq a+1$, and $\alt^r_\sigma(\HH)\geq a+b$ for $r\geq 4$, 
which clearly complete the proof.  
Consider an arbitrary bijection $\sigma:[n]\longrightarrow [n]$. 
Setting $X=(x_1,\ldots,x_n)$ with 
$$x_i=\left\{\begin{array}{ll}
\omega^i     &  i\in\{1,\ldots,r(k-1)\}\\
0		&  \mbox{otherwise},  
\end{array}\right.$$
concludes $\alt(X)=r(k-1)$ and $|X^i|=k-1$ for each $i\in[r]$ implying $E(\HH[\sigma(X^i)])=\varnothing$ for each $i\in[r]$. 
Therefore, in view of the definition of $\alt^r_\sigma(\HH)$, we have $\alt^r_\sigma(\HH)\geq r(k-1)$. 
Since $\sigma$ is chosen arbitrarily, this implies $\alt^r(\HH)\geq r(k-1)$.  Note that if $a\leq k-1$, then $\HH=([n],{[n]\choose k})$.  
Now, for an arbitrary bijection $\sigma:[n]\longrightarrow [n]$, if $\alt(X)\geq r(k-1)+1$, then for at least one $j$, we must have $|X^j|\geq k$
which implies $E(\HH[\sigma(X^i)])\neq\varnothing$. This concludes $\alt^r(\HH)=\alt^r([n],{[n]\choose k})=r(k-1)$ completing the proof in this case. 
Henceforth, we suppose that $a\geq k$. Set $A=\{n-a+1,\ldots,n\}$. 
Let $\sigma:[n]\longrightarrow [n]$ be an arbitrary bijection. 
Also, let $\sigma^{-1}(A)=\{j_1,\ldots,j_a\}$,   
where $1\leq j_1<\cdots<j_a\leq n$. 
The following three different cases will be distinguished. 
\begin{itemize}
\item {\bf Case $\mathbf{r=2}:$} 
Define $X=(x_1,\ldots,x_n)$ such that 
$$x_i=\left\{
\begin{array}{ll}
\omega     &  i=j_q\in\sigma^{-1}(A)\mbox{ and } q \mbox{ is even}\\
\omega^2 &  i=j_q\in\sigma^{-1}(A)\mbox{ and } q \mbox{ is odd}\\
0		&  \mbox{otherwise}.
\end{array}\right.
$$
Accordingly, we have $\alt(X)=|A|=a$ and $\sigma(X^i)\subseteq A$ for each $i\in[r]$. 
This concludes $E(\HH[\sigma(X^i)])=\varnothing$ for each $i\in[r]$ 
implying that $\alt^r_\sigma(\HH)\geq a$. Therefore, since $\sigma$ is arbitrary, we have 
$\alt^r(\HH)\geq a$. Note that we already proved that $\alt^r(\HH)\geq r(k-1)$ which completes the proof for $r=2$. 

\item{\bf Case $\mathbf{r=3:}$}
Consider a fixed $i_0\in [n]\setminus A$. 
Set $X=(x_1,\ldots,x_n)$ such that 
$$x_i=\left\{
\begin{array}{ll}
\omega     &  i=j_q\in\sigma^{-1}(A)\mbox{ and } q \mbox{ is even}\\
\omega^2 &  i=j_q\in\sigma^{-1}(A)\mbox{ and } q \mbox{ is odd}\\
\omega^3 &  i=\sigma^{-1}(i_0)\\
0		&  \mbox{otherwise.}
\end{array}\right.$$
Note $\alt(X)=a+1$ and moreover, $\sigma(X^1), \sigma(X^2)\subseteq A$ and $|\sigma(X^3)|=1$ which clearly implies 
$E(\HH[\sigma(X^i)])=\varnothing$ for each $i\in[r]$. 
Therefore, in view of the definition of $\alt^r_\sigma(\HH)$, we have  $\alt^r_\sigma(\HH)\geq \alt(X)=a+1$. Since $\sigma$ is chosen arbitrarily, we have 
$\alt^r(\HH)\geq a+1$. Note that we have already proved that $\alt^r(\HH)\geq r(k-1)$, which completes the proof for $r=3$. 

\item{\bf Case $\mathbf{r\geq 4:}$}
Consider a fixed $B=\{1,\ldots,b\}\subseteq [n]\setminus A$ such that 
$b=\min\left\{n-a,(r-2)(k-1)\right\}$.  
Also, let $\sigma^{-1}(B)=\{l_1,\ldots,l_b\}$, where $1\leq l_1<\cdots<l_b\leq n$. 
Define $X=(x_1,\ldots,x_n)$ such that 
$$x_i=\left\{
\begin{array}{ll}
\omega      &  i=j_q\in\sigma^{-1}(A)\mbox{ and } q \mbox{ is even}\\
\omega^2  &  i=j_q\in\sigma^{-1}(A)\mbox{ and } q \mbox{ is odd}\\
\omega^{3+t}  &  i=l_q\in\sigma^{-1}(B)\mbox{ and } q\equiv t\in\{0,\ldots,r-3\}\; (\mod r-2)\\
0		 &  \mbox{otherwise}.
\end{array}\right.
$$
Clearly, $\alt(X)=|A|+b$ and $\sigma(X^1), \sigma(X^2)\subseteq A$ and $|\sigma(X^i)|\leq k-1$ for each $i\geq 3$ which  conclude 
$E(\HH[\sigma(X^i)])=\varnothing$ for each $i\in[r]$.   
Therefore, in view of the definition of $\alt^r_\sigma(\HH)$, we have $\alt^r_\sigma(\HH)\geq a+b$. Since $\sigma$ is chosen arbitrarily, we have 
$\alt^r(\HH)\geq \alt(X)=a+b$. This result in combination with the fact that $\alt^r(\HH)\geq r(k-1)$ implies the proof when $r\geq 4$. 
\end{itemize}
}\end{proof}

Suppose that $r\geq 4$, $a\geq rk-1$, and  $n\geq (r-1)(k-1)+a$. For simplicity of notation, set $\HH=\HH(n,k,a)$. 
On the one hand, 
Observation~\ref{ecrH(n,k,a)} for $\ss=(1,\ldots,1)$ yields $\ecd^r(\HH)=n-a$ and $\cd^r(\HH)=n-(r-1)(k-1)-a$, 
but on the other hand,  Lemma ~\ref{altcompare} implies $n-\alt^r(\HH) \leq n- a-(r-2)(k-1)$. Therefore,  
$\ecd^r(\HH)-\cd^r(\HH)=(r-1)(k-1)$ and 
$\ecd^r(\HH)-(n -\alt^r(\HH))\geq {(r-2)(k-1)},$
which clearly implies the following observation.
\begin{observation}
For $r\geq 4$, $a\geq rk-1$, and  $n\geq (r-1)(k-1)+a$, the values of ${1\over r-1}\big(\ecd^r(\HH)-\cd^r(\HH)\big)$ and ${1\over r-1}\big(\ecd^r(\HH)-(n -\alt^r(\HH))\big)$
can be made as large as desired by making $k$ large enough. 
\end{observation} 

For the hypergraph $\HH=\HH(n, k , a)$, the following proposition demonstrates that 
$\left\lceil{\ecd^r(\HH)\over r-1}\right\rceil$ provides an exact lower bound for the 
chromatic number of $\KG^r(n,k,a)$ in some cases, while, 
in view of the previous observation,  neither 
$\left\lceil{\cd^r(\HH)\over r-1}\right\rceil$ nor 
$\left\lceil{n-\alt^r(\HH)\over r-1}\right\rceil$ do that. 
Note that the following proposition can be considered as a generalization of the 
Alon-Frankl-Lov{\'a}sz theorem as well. 

\begin{proposition} \label{proposition1}
Let $n, r,k,$ and $a$ be  positive integers, where $r,k\geq 2$, $n> a$, and $n \geq rk$. 
Then we have 
$$\chi(\KG^r(n,k,a))=\left\lceil \dfrac{n-\max\{r(k-1), a\}}{r - 1}\right\rceil,$$ 
provided that either $a\leq 2(k-1)$ or $a\geq rk-1$. 
Moreover, for $2k -1 \leq a\leq rk-2$, we have 
$$\left\lceil\dfrac{n - r(k-1) - \lfloor{a\over k}\rfloor}{r-1}\right\rceil\leq \chi(\KG^r(n,k,a))\leq \left\lceil\dfrac{n - \max\{r(k-1),a\}}{r-1}\right\rceil.$$ 
\end{proposition}
\begin{proof}{
The cases $a<k$ and $a\geq rk-1$ have been already proved by Corollary~\ref{corchrom} when $\ss=(1,\ldots,1)$. 
Also, by Lemma~\ref{upper} and Observation~\ref{ecrH(n,k,a)}, we have  
$$\left\lceil\dfrac{n - r(k-1) - \lfloor{a\over k}\rfloor}{r-1}\right\rceil\leq \chi(\KG^r(n,k,a))\leq \left\lceil\dfrac{n - \max\{r(k-1), a\}}{r-1}\right\rceil,$$  
provided that $k \leq a\leq rk-2$. Therefore, to complete the proof, we must show that 
$$\chi(\KG^r(n,k,a))\geq\left\lceil \dfrac{n-r(k-1)}{r - 1}\right\rceil,$$
provided that $k\leq a\leq 2(k-1)$. 
Since $\KG^r(n,k,a)=\KG^r(\HH(n,k,a))$, in view of Inequality~\ref{alishahi&haji inequality}, we have 
$$\chi(\KG^r(n,k,a)) \geq \left\lceil\dfrac{|V(\HH(n,k,a))| - \alt^r_I(\HH(n,k,a))}{r-1}\right\rceil,$$
where $I:[n]\longrightarrow[n]$ is the identity bijection. Hence, we have the proof concluded if we prove that $\alt^r_I(\HH(n,k,a))\leq r(k-1)$. 
To compute $\alt^r_I(\HH(n,k,a))$, we must find the maximum possible value of $\alt(X)$ over all $X\in(Z_r\cup\{0\})^n$ such that none of $X^i$'s contains an edge of $\HH(n,k,a)$. 
It is clear that this maximum also happens among $X$'s with $\alt(X)=|X|$, that is, 
$$\alt^r_I(\HH)=\max\left\{\alt(X)\colon X\in ( Z_r\cup\{0\})^n,\ |X|=\alt(X), \mbox{ and }  
E(\HH[X^i])=\varnothing\mbox{ for each } i\in[r]\right\}.$$
Note $\alt(X)=|X|$ implies that  none of  $X^i$'s contains consecutive elements of $[n]$ and therefore, $|X^i\cap\{n-a+1,\ldots,n\}|\leq \lceil{a\over 2}\rceil\leq k-1$ for each $i$. 
Using this observation along with the fact that if $\alt(X)=|X|\geq r(k-1)+1$, then there is at least one $i\in [r]$ for which $|X^i|\geq k$ simply concludes   
$\alt^r_I(\HH(n,k,a))\leq r(k-1)$ as desired. 
} \end{proof}

We here pose the following conjecture asserting that  the hypothesis $a\leq 2(k-1)$ or $a\geq rk-1$ is superfluous in the preceding proposition.  
\begin{conjecture}\label{conjecture}
Let $n, r,k,$ and $a$ be  positive integers, where $r,k\geq 2$, $n> a$, and $n \geq rk$. We have 
$$\chi(\KG^r(n,k,a))=\left\lceil \dfrac{n-\max\{r(k-1), a\}}{r - 1}\right\rceil.$$  
\end{conjecture}
This conjecture is strongly supported by Proposition~\ref{proposition1}.
To be more specific, note that the validity of this conjecture is already verified for $a\leq 2(k-1)$,  $a\geq rk-1$, and, in particular, for 
$r=2$ by Proposition~\ref{proposition1}. Therefore, it is left open just for $2k-1\leq a\leq rk-2$ with $r\geq 3$. 
Even more, if $a\leq r(k-1)$ and $r$ is a power of $2$, then  
we can deduce Conjecture~\ref{conjecture} from
a result by Alon, Drewnowski, and {\L}uczak~\cite{MR2448565}.  
For $t\geq 2$, a subset  $A$ of $[n]$  is called $t$ stable if $t\leq |x-y|\leq n-t$ 
for each $x\neq y\in A$. 
The induced subhypergraph of $\KG^r(n,k)$ by the set of all $t$-stable vertices is called 
the {\it $t$-stable Kneser hypergraph $\KG^r(n,k)_{t-stab}$}. 
Alon, Drewnowski, and {\L}uczak~\cite{MR2448565} proved that
$\chi(\KG^r(n,k)_{r-stab})=\chi(\KG^r(n,k))$ provided that $r$ is a power of $2$. They also 
conjectured that this result is true for general $r$. One can simply see that the hypergraph 
$\KG^r(n,k,a)$ in the statement of Conjecture~\ref{conjecture} is a super-hypergraph 
of $\KG^r(n,k)_{r-stab}$ provided that $a\leq r(k-1)$, which implies that Conjecture~\ref{conjecture} 
is true provided that  $a\leq r(k-1)$ and $r$ is a power of $2$ . 
In this point of view, we can consider this conjecture as a weak version of the Alon-Drewnowski-{\L}uczak conjecture.

For a hypergraph $\HH$, though the lower bound of $\chi(\KG^r(\HH))$ which is based on $\ecd^r(\HH)$ might be 
appreciably better than the lower bound based on $\alt^r(\HH)$, in general we cannot decide which one is better. 
There are several examples in which the lower bound of $\chi(\KG^r(\HH))$ based on $\alt^r(\HH)$ is much better than the lower bound gained by $\ecd^r(\HH)$.
As an instance, for positive integers $n$, $r$, and $k$, 
let $K =K_{n,...,n}$ be a graph with vertex set $[rn]$ and edge set $E(K) = \left\{e \in {[rn]\choose 2}: e \nsubseteq \{1 + (t-1)n,...,tn\} \quad \forall t \in [r]\right\}$. 
By the definition of $K$, it is observed that $V(K)$ 
is partitioned into $r$ pairwise disjoint sets of size $n$ such that none of them contains any edge 
of $K$. So, $\ecd^r(K) = 0$, while  it can be verified that 
$$rn-\alt^r(K)\geq rn-\alt^r_I(K)=\left\{\begin{array}{ll}
{r\over 2}n & \mbox{ $r$ is even}\\ \\
\lceil{r\over 2}\rceil n-1 & \mbox{ $r$ is odd,}
\end{array}\right.$$
where $I:[rn]\longrightarrow [rn]$ is the identity bijection.  
\begin{proposition}
Let $r$ be a fixed integer where $r\geq 2$. 
For a given hypergraph $\mathcal{F}$, 
it is notable that though $\left\lceil{\ecd^r(\mathcal{F})\over r-1}\right\rceil$ 
provides a lower bound for the chromatic number of  $\KG^r(\mathcal{F})$, 
computing $\ecd^r(\mathcal{F})$ is an $NP$-hard problem. 
\end{proposition}
The proof of this observation is almost the same as that of  Proposition~6 in~\cite{Meunier2014}.
However, for the sake of completeness,  we sketch the proof in the following. 
  
Let $r$ be a fixed integer, where $r\geq 2$. Also, let $G$ be a given graph for
which we want to compute the independence number $\alpha(G)$. 
Define $\bar{G}$ to be a graph which is constructed    
by the union of $r$ vertex-disjoint copies of $G$ such that any two vertices from different 
copies are adjacent. It is~not difficult to see that 
$$\ecd^r(\bar{G})=r\big(|V(G)|-\alpha(G)\big).$$
This observation implies that computing the equitable colorability defect of 
graphs is at least as difficult as 
computing the independence number of graphs. Since computing the independence number of graphs is 
known as an NP-hard problem, we have the proof completed. 
 
\begin{remark}
After this paper became available online, using different methods, Aslam {\it et al.}~\cite{2017arXiv171203456A} proved Conjecture~\ref{conjecture} in some cases.  
\end{remark}

\section*{Acknowledgments}  
We would like to thank professor Fr\'ed\'eric Meunier for his comments that helped  
to improve the presentation of the paper. 
Also, we thank the two anonymous reviewers for their careful reading of our manuscript and their constructive comments. 

%%%%%%%%%%%%%%%%%%%%%%%%%%%%%%%%%%%%%%%%%%%%%
%%%%%%%%%%%%%%%%%%%%%%%%%%%%%%%%%%%%%%%%%%%%%
%\bibliographystyle{plain}
%\bibliography{MyReferences}

\end{document}